\documentclass[11pt]{article}
\usepackage[matrix,arrow,curve,cmtip]{xy}
\usepackage{graphicx}
\usepackage{amssymb}
\usepackage{latexsym}
\usepackage{theorem}
\usepackage[a4paper,width=140mm,height=220mm]{geometry}
\usepackage{titlesec}
\titleformat{\section}[hang]%
{\bf\filcenter\large}{\thesection.}{1ex}{}%
\titleformat{\subsection}[runin]%
{\bfseries\normalsize}{}{0ex}{}

\def\to{\mbox{$\xymatrix@1@C=5mm{\ar@{->}[r]&}$}}
\def\tto{\mbox{$\xymatrix@1@C=5mm{\ar@{=>}[r]&}$}}
\def\biar{\mbox{$\xymatrix@1@C=5mm{\ar@<1.5mm>[r]\ar@<-0.5mm>[r]&}$}}
\def\adjar{\mbox{$\xymatrix@1@C=5mm{\ar@<1.5mm>@{<-}[r]\ar@<-0.5mm>[r]&}$}}
\def\iso{\mbox{$\xymatrix@1@C=6mm{\ar@{->}[r]^{\sim}&}$}}
\def\distsign{\begin{picture}(0,0)\put(0,0){\circle{4}}\end{picture}}
\def\dist{\mbox{$\xymatrix@1@C=5mm{\ar@{->}[r]|{\distsign}&}$}}

\newtheorem{theorem}{Theorem}[section]
\newtheorem{lemma}[theorem]{Lemma}

{\theorembodyfont{\upshape}}
{\theorembodyfont{\upshape}}
\newcommand{\proof}{\noindent {\it Proof\ }: }

\def\endofproof{$\mbox{ }\hfill\Box$\par\vspace{1.8mm}\noindent}

\def\:{\colon}
\def\1{\mathbf{1}}

\def\2{\mathbf{2}}
\def\inv{^{-1}}

\def\Matr{{\sf Matr}}
\def\Dist{{\sf Dist}}
\def\Cat{{\sf Cat}}

\def\Q{\mathcal{Q}}

\def\C{\mathcal{C}}

\def\V{\mathcal{V}}

\def\lim{\mathop{\rm lim}}
\def\bbA{\mathbb{A}}
\def\bbB{\mathbb{B}}
\def\bbC{\mathbb{C}}

\def\bbP{\mathbb{P}}

\def\bbX{\mathbb{X}}
\def\bbY{\mathbb{Y}}
\def\bbT{\mathbb{T}}
\def\tensor{\otimes}

\def\<{\langle}
\def\>{\rangle}
\def\eqref#1{(\ref{#1})}
\def\defeq{:=}
\def\name#1{[#1]}
\def\pb{%
\setlength{\unitlength}{1ex}
\begin{picture}(1,1.5)(1,0.2)
\put(0,0){\line(1,0){2}} \put(2,0){\line(0,1){2}}
\end{picture}}
\def\otspam{%
\setlength{\unitlength}{1ex}
\begin{picture}(3.7,1)(0,0)
\put(0.6,0.3){\mbox{\rotatebox[origin=c]{180}{$\mapsto$}}}
\end{picture}}

\def\slice#1{%
\setlength{\unitlength}{1ex}
\begin{picture}(2.5,1)(0,0)
\put(-0.2,-0.5){\mbox{$/$}}
\put(1,-1){\mbox{${#1}$}}
\end{picture}}

\title{Exponentiable functors \\ between quantaloid-enriched categories}
\author{Maria Manuel Clementino\footnote{%
     Centro de Matem\'atica, 
     Universidade de Coimbra, Portugal. Email: 
     {\tt mmc@mat.uc.pt}}, 
\ Dirk Hofmann\footnote{%
     Departamento de Matem\'atica, 
     Universidade de Aveiro, Portugal. Email: 
     {\tt dirk@mat.ua.pt}} 
\ and Isar Stubbe\footnote{%
     Centro de Matem\'atica, 
     Universidade de Coimbra, Portugal. Email: 
     {\tt isar@mat.uc.pt}}}
\date{August 17, 2006}

\begin{document}

\maketitle

\begin{quote}{\small
{\bf Abstract.} Exponentiable functors between quantaloid-enriched 
categories are characterized in elementary terms. The proof goes as 
follows: the elementary conditions on a given functor translate into existence 
statements for certain adjoints that obey some lax commutativity; 
this, in turn, is precisely what is needed to prove the existence of 
partial products with that functor; so that the functor's exponentiability follows from the works of Niefield [1980] and Dyckhoff and Tholen [1987].
\\
{\bf Keywords:} quantaloid, enriched category, exponentiability, 
partial product. 
\\
{\bf MSC 2000 Classification:} 06F07, 18A22, 18D05, 18D20 
}\end{quote}

\section{Introduction}\label{A}
The study of exponentiable morphisms in a category $\cal C$, in 
particular of exponentiable functors between (small) categories 
(i.e.\ Conduch\'e fibrations), has a long history; see [Niefield, 
2001] for a short account. Recently M. M. Clementino and D. Hofmann 
[2006] found simple necessary-and-sufficient conditions for the 
exponentiability of a functor between $\cal V$-enriched categories, 
where $\cal V$ is a symmetric quantale which has its top element as 
unit for its multiplication and whose underlying sup-lattice is a 
locale. Our aim here is to prove the following 
characterization of the exponentiable functors between $\Q$-enriched 
categories, where now $\Q$ is {\it any (small) quantaloid}, thus 
considerably generalizing the aforementioned result of [Clementino 
and Hofmann, 2006].
\begin{theorem}\label{1}
A functor $F\:\bbA\to\bbB$ between $\Q$-enriched categories is 
exponentiable, i.e.\ the functor ``product with $F$''
$$-\times F\:\Cat(\Q)\slice{\bbB}\to\Cat(\Q)\slice{\bbB}$$
admits a right adjoint, if and only if the following two conditions hold:
\begin{enumerate}
\item\label{e1} for every $a,a'\in\bbA$ and $\bigvee_if_i\leq\bbB(Fa',Fa)$, $$\Big(\bigvee_if_i\Big)\wedge\bbA(a',a)=\bigvee_i\Big(f_i\wedge\bbA(a',a)\Big),$$
\item\label{e2} for every $a,a''\in\bbA$, $b'\in\bbB$, $f\leq\bbB(b',Fa)$ and $g\leq\bbB(Fa'',b')$,
$$(g\circ f)\wedge\bbA(a'',a)=\bigvee_{a'\in F\inv b'}\Big((g\wedge\bbA(a'',a'))\circ(f\wedge\bbA(a',a))\Big).$$
\end{enumerate}
\end{theorem}
These conditions are ``elementary'' in the sense that they are 
simply equalities (of infima, suprema and compositions) of morphisms 
in the base quantaloid $\Q$. The second condition is precisely what 
[Clementino and Hofmann, 2006] had too, albeit in their more 
restrictive setting; but they did not discover the first condition 
{\em an sich}: because it is obviously always true if the base 
category is a locale.
\par
The proof of our theorem goes as follows. In section \ref{B} we 
first translate conditions \ref{1}--\ref{e1} and \ref{1}--\ref{e2} into existence 
statements for certain adjoints obeying some lax commutativity. 
Next, in section \ref{C}, we show that these latter adjoints are 
precisely what is needed to prove the existence of partial products 
in $\Cat(\Q)$ over $F\:\bbA\to\bbB$. The result then follows from R. 
Dyckhoff and W. Tholen's [1987] observation, complementary to S. 
Niefield's [1982] work, that a morphism $f\:A\to B$ in a category 
$\cal C$ with finite limits is exponentiable if and only if $\C$ 
admits partial products over $f$.
\\[2mm]
{\bf Acknowledgement.} This work was done when Isar Stubbe was a 
post-doctoral researcher at the Centre for Mathematics of the University 
of Coimbra.

\section{Preliminaries}\label{A2}
For the basics on $\Q$-enriched categories we refer to [Stubbe, 
2005]; all our notations are as in that paper. Here we shall just observe that $\Cat(\Q)$ has pullbacks and
a terminal -- and therefore all finite limits [Borceux, 1994, 
Proposition 2.8.2] -- and fix some notations.
\par
The terminal object in $\Cat(\Q)$, write it as $\bbT$, has:
\begin{itemize}
\item objects: $\bbT_0=\Q_0$, with types $tX=X$,
\item hom-arrows: $\bbT(Y,X)=\top_{X,Y}=$ the top element of $\Q(X,Y)$.
\end{itemize}
For two functors $F\:\bbA\to\bbC$ and $G\:\bbB\to\bbC$ with common 
codomain, their pullback $\bbA\times_{\bbC}\bbB$ has:
\begin{itemize}
\item objects: $(\bbA\times_{\bbC}\bbB)_0=\{(a,b)\in\bbA_0\times\bbB_0)\mid Fa=Gb\}$ with $t(a,b)=ta=tb$,
\item hom-arrows: $(\bbA\times_{\bbC}\bbB)((a',b'),(a,b))=\bbA(a',a)\wedge\bbB(b',b)$,
\end{itemize}
and comes with the obvious projections. All verifications are 
entirely straightforward.
\par
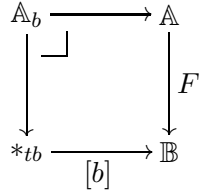
\begin{figure}
$$\xy\xymatrix@=8ex{
\bbA_b\ar[r]\ar[d] & \bbA\ar[d]^F \\
{*_{tb}}\ar[r]_{\name{b}} & \bbB}\POS(4.5,-4)\drop{\pb}\endxy$$ 
\caption{a specific pullback}\label{f0}
\end{figure}

For an $X\in\Q$, the one-object $\Q$-category with hom-arrow $1_X$ 
is written as $*_X$. There is an obvious bijection between the 
objects of type $X$ in some $\Q$-category $\bbB$ and the functors 
from $*_X$ to $\bbB$. Thus, let $\name{b}\:*_{tb}\to\bbB$ stand for 
the functor ``pointing at'' $b\in\bbB$. Given a functor 
$F\:\bbA\to\bbB$ and an object $b\in\bbB$ in its codomain, we shall 
write $\bbA_b$ for the pullback in figure \ref{f0}.
That is to say, $\bbA_b$ has
\begin{itemize}
\item objects: $(\bbA_b)_0=F\inv b=\{a\in\bbA\mid b=Fa\}$, all of type $tb$,
\item hom-arrows: $\bbA_b(a',a)=1_{tb}\wedge\bbA(a',a)$.
\end{itemize}
Note that $\bbA_b=\emptyset$ if and only if $b\not\in F(\bbA)$. 

\section{Adjoints obeying a lax commutativity}\label{B}
In this section we shall translate the elementary conditions in \ref{1} 
into existence statements of certain adjoints obeying some lax 
commutative diagrams.
\par
\begin{lemma}\label{2}
For a functor $F\:\bbA\to\bbB$ between $\Q$-categories, the 
following are equivalent conditions:
\begin{enumerate}
\item condition \ref{1}--\ref{e1} holds,
\item for every $a,a'\in\bbA$, the order-preserving map 
\begin{equation}\label{3}
\downarrow\bbB(Fa',Fa)\to\Q(ta,ta')\:f\mapsto f\wedge\bbA(a',a)
\end{equation}
has a right adjoint,
\item\label{2.a} for every $b,b'\in F(\bbA)$, the order-preserving map 
\begin{equation}\label{4}
\downarrow\bbB(b',b)\to\Matr(\Q)(\bbA_b,\bbA_{b'})\:f\mapsto\Big(f\wedge\bbA(a',a)\Big)_{(a,a')\in\bbA_b\times\bbA_{b'}}
\end{equation}
has a right adjoint.
\item\label{2.b} for every $b,b'\in F(\bbA)$, the order-preserving map 
\begin{equation}\label{4.1}
\downarrow\bbB(b',b)\to\Dist(\Q)(\bbA_b,\bbA_{b'})\:f\mapsto 
\Big(f\wedge\bbA(a',a)\Big)_{(a,a')\in\bbA_b\times\bbA_{b'}}
\end{equation}
has a right adjoint.
\item\label{a1} for every $b,b'\in\bbB$, the order-preserving map in \eqref{4.1}
has a right adjoint.
\end{enumerate}
\end{lemma}
\proof The equivalence of the first two statements is trivial: an 
order-preserving map between complete lattices has a right adjoint 
if and only if it preserves arbitrary suprema. 
\par
Next, if we use $g\mapsto g^F$ as generic notation for the right 
adjoints to the maps in \eqref{3}, then
$$M\mapsto M^F\defeq\bigwedge\{M(a',a)^F\mid (a,a')\in\bbA_b\times\bbA_{b'}\}$$
is the right adjoint to the map in \eqref{4}. Conversely, if $M\mapsto 
M^F$ is the right adjoint to the map in \eqref{4}, then for any 
$a,a'\in\bbA$
$$g\mapsto g^F\defeq\left(T^{(a,a')}(g)\right)^F$$
is the right adjoint to the map in \eqref{3}, with $T^{(a,a')}(g)$ 
standing for the $\Q$-matrix from $\bbA_{Fa}$ to $\bbA_{Fa'}$ all of 
whose elements are set to the top element in $\Q(ta,ta')$
except for the element indexed by $(a,a')$ which is set to $g$. 
\par
The equivalence of \ref{2.a} and \ref{2.b} follows straightforwardly from two 
facts: First, the {\it matrix}
$$\widehat{f}\defeq\Big(f\wedge\bbA(a',a)\Big)_{a\in\bbA_b,a'\in\bbA_{b'}}$$
is always a {\it distributor} from $\bbA_b$ to $\bbA_{b'}$: because 
for any $a,a_1\in\bbA_b$ and $a',a_1'\in\bbA_{b'}$ it is automatic 
that 
\begin{eqnarray*} 
\widehat{f}(a',a_1)\circ\bbA_b(a_1,a)
 & = & \Big(f\wedge\bbA(a',a_1)\Big)\circ\Big(1_{ta}\wedge\bbA(a_1,a)\Big) \\
 & \leq & \Big(f\circ 1_{ta}\Big)\wedge\Big(\bbA(a',a_1)\circ\bbA(a_1,a)\Big) \\
 & \leq & f\wedge\bbA(a',a) \\
 & = & \widehat{f}(a',a)
\end{eqnarray*}
and similarly 
$\bbA_{b'}(a'_1,a')\circ\widehat{f}(a',a)\leq\widehat{f}(a',a)$. And 
second, the inclusion 
$$\Dist(\Q)(\bbA_{b'},\bbA_b)\to\Matr(\Q)(\bbA_{b'},\bbA_b)\:\Phi\mapsto\Phi$$
has both a left and a right adjoint; namely, its left adjoint is 
$M\mapsto \bbA_{b'}\circ M\circ\bbA_b$ and its right adjoint is
$M\mapsto [\bbA_{b'},\{\bbA_b,M\}]$. (In both expressions, $\bbA_{b'}$ and 
$\bbA_b$ are viewed as monads in the quantaloid $\Matr(\Q)$, and we compute  
composition, resp.\ lifting and extension, of matrices.) Hence both triangles in figure 
\ref{f6} commute and both solid arrows are left adjoints, so it 
follows that one dashed arrow is a left adjoint if and only if the 
other one is.
\begin{figure}
$$\xymatrix@R=10ex@C=2ex{ & \downarrow\bbB(b',b)\ar@{.>}[dl]_{f\mapsto\widehat{f}}\ar@{.>}[dr]^{f\mapsto\widehat{f}} \\
\Matr(\Q)(\bbA_b,\bbA_{b'})\ar@<-1mm>[rr]_{M\mapsto \bbA_{b'}\circ 
M\circ \bbA_b} & & 
\Dist(\Q)(\bbA_b,\bbA_{b'})\ar@<-1mm>[ll]_{\Phi\otspam\Phi}}$$ 
\caption{a diagram for the proof of \ref{2}}\label{f6}
\end{figure}
\par
Finally, the only difference between the fourth and the fifth 
statement is that in the latter it may be that $\bbA_b$ or 
$\bbA_{b'}$ is empty; but then $\Dist(\Q)(\bbA_b,\bbA_{b'})$ is a 
singleton (containing the empty distributor) in which case the right 
adjoint to \eqref{4.1} always exists.
\endofproof
\par
In the statement of the next lemma we shall write
\begin{equation}\label{4.1b}
\downarrow\bbB(b',b)
\xymatrix@C=8ex{\ar@{}[r]|{\perp}\ar@<1mm>@/^2mm/[r]^{f\mapsto\widehat{f}}
 & \ar@<1mm>@/^2mm/[l]^{\Phi^F\otspam\Phi}} \Dist(\Q)(\bbA_b,\bbA_{b'})
\end{equation}
for the adjunctions (one for each pair $(b,b')$ of objects of $\bbB$) that \ref{2}--\ref{a1} alludes to.
\begin{lemma}\label{5}
For a functor $F\:\bbA\to\bbB$ between $\Q$-categories for which the 
equivalent conditions in \ref{2} hold, the following are equivalent 
conditions:
\begin{enumerate}
\item condition \ref{1}--\ref{e2} holds,
\item for every $a,a''\in\bbA$ and $b'\in\bbB$, the diagram 
in figure \ref{f1}, in which the horizontal arrows are given by 
composition (in $\Dist(\Q)$, resp.\ $\Q$), the left vertical arrow is
\begin{equation}\label{5.1}
(f,g)\mapsto\Big(\left(f\wedge\bbA(a',a)\right)_{a'\in\bbA_{b'}},\left(g\wedge\bbA(a'',a')\right)_{a'\in\bbA_{b'}}\Big)
\end{equation}
and the right vertical arrow is as in \eqref{3}, is lax 
commutative as indicated,
\begin{figure}[t]
$$\xy\xymatrix@u@C=5ex@R=0ex{
\downarrow\bbB(b',Fa)\times\downarrow\bbB(Fa'',b')\ar[rr]\ar[d]_(0.58){-\circ-} & & \Dist(\Q)(*_{ta},\bbA_{b'})\times\Dist(\Q)(\bbA_{b'},*_{ta''})\ar[d]^(0.6){-\tensor-} \\
\downarrow\bbB(Fa'',Fa)\ar[r] & \Q(ta,ta'')\ar@{=}[r]& 
\Dist(\Q)(*_{ta},*_{ta''})} \POS(30,14)\drop{\geq}\endxy$$ 
\caption{the diagram for \ref{5}--2}\label{f1}
\end{figure}
\item for every $b,b''\in F(\bbA)$ and $b'\in\bbB$, the diagram in figure \ref{f2} is lax commutative as indicated,
\begin{figure}[t]
$$\xy\xymatrix@u@C=10ex@R=1ex{
\downarrow\bbB(b',b)\times\downarrow\bbB(b'',b')\ar[r]^{(\widehat{-})\times(\widehat{-})}\ar[d]_(0.58){-\circ-} & \Dist(\Q)(\bbA_b,\bbA_{b'})\times\Dist(\Q)(\bbA_{b'},\bbA_{b''})\ar[d]^(0.6){-\tensor-} \\
\downarrow\bbB(b'',b)\ar[r]_{(\widehat{-})} & \Dist(\Q)(\bbA_b,\bbA_{b''})} 
\POS(35,10)\drop{\geq}\endxy$$ \caption{the diagram for \ref{5}--3 
and \ref{5}--4}\label{f2}
\end{figure}
\item for every $b,b',b''\in\bbB$, the diagram in figure \ref{f2} is lax commutative as indicated,
\item\label{a2} for every $b,b',b''\in\bbB$, the diagram in figure \ref{f3} is 
lax commutative as indicated.
\begin{figure}[t]
$$\xy\xymatrix@R=10ex@C=11ex{
\Dist(\Q)(\bbA_b,\bbA_{b'})\times\Dist(\Q)(\bbA_{b'},\bbA_{b''})\ar[d]_{(-)^F\times(-)^F}\ar[r]^(0.6){-\tensor-} & \Dist(\Q)(\bbA_b,\bbA_{b''})\ar[d]^{(-)^F} \\
\downarrow\bbB(b',b)\times\downarrow\bbB(b'',b')\ar[r]_(0.58){-\circ-} & 
\downarrow\bbB(b'',b)} \POS(35,-12)\drop{\leq}\endxy$$ \caption{the 
diagram for \ref{5}--\ref{a2}}\label{f3}
\end{figure}
\end{enumerate}
\end{lemma}
\proof First it is easily verified, in an analogous manner as in the previous proof, that the map in \eqref{5.1} is well-defined, i.e.\ that we indeed defined {\it distributors}
$$\left(f\wedge\bbA(a',a)\right)_{a'\in\bbA_{b'}}\mbox{, resp.\ }\left(g\wedge\bbA(a'',a')\right)_{a'\in\bbA_{b'}}$$
from $*_{ta}$ to $\bbA_{b'}$, resp.\ from $\bbA_{b'}$ to $*_{ta''}$. Now the equivalence of the first two statements is immediate; the 
``oplax commutativity'' of the diagram in figure \ref{f1} is always 
true, thus explaining why in \ref{1}--\ref{e2} there is an equality instead of an 
inequality. That the second and the third statement are equivalent, 
is because all order-theoretic operations on a distributor are done ``elementwise''; and 
the third and fourth are equivalent because in case $\bbA_b$ or 
$\bbA_{b''}$ is empty, $\Dist(\Q)(\bbA_b,\bbA_{b''})$ is a 
singleton, hence all is trivial. Finally, the equivalence of the two 
last statements follows from the respective vertical arrows being 
adjoint.
\endofproof

\section{Partial products}\label{C}
In this section we link the conditions in \ref{2} and \ref{5} on a functor $F\:\bbA\to\bbB$ to the existence of so-called partial products in $\Cat(\Q)$ with $F$: this completes the proof of \ref{1}.
\par
\begin{figure}[t]
$$\xymatrix@=5ex{
\bbC & & \bbP\times_{\bbB}\bbA\ar[ll]_E\ar[rr]\ar[dd] & & \bbA\ar[dd]^F \\
  & \bbP'\times_{\bbB}\bbA\ar@<-1mm>[urrr]\ar[dd]\ar[ul]^{E'}\ar@{.>}[ur]^{K\times_{\bbB}1_{\bbA}} \\
  & & \bbP\ar[rr]^P & & \bbB \\
  & \bbP'\ar@<-1mm>[urrr]_{P'}\ar@{.>}[ur]^K}$$
\caption{the definition of a partial product}\label{f4} 
\end{figure}
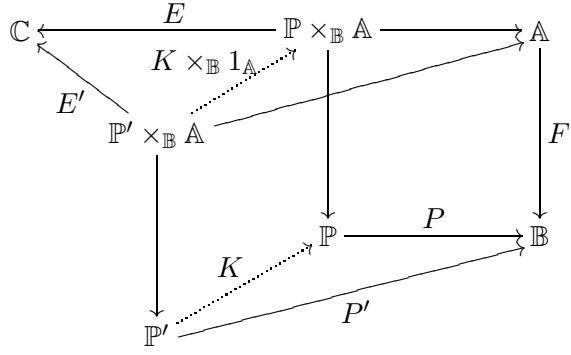

First recall R. Dyckhoff and W. Tholen's [1987] definition (which they 
gave for any morphism $f\:A\to B$ and any object $C$ in any category 
$\C$ with finite limits, but here it is for $\Q$-cat\-e\-gor\-ies): 
the {\em partial product} of a functor $F\:\bbA\to\bbB$ with a 
$\Q$-category $\bbC$ is a $\Q$-category $\bbP$ together with 
functors $P\:\bbP\to\bbB$, $E\:\bbP\times_{\bbB}\bbA\to\bbC$ such 
that, for any other $\Q$-category $\bbP'$ and functors 
$P'\:\bbP'\to\bbB$, $E'\:\bbP'\times_{\bbB}\bbA\to\bbC$ there exists 
a unique functor $K\:\bbP'\to\bbP$ satisfying $P\circ K=P'$ and 
$E\circ(K\times_{\bbB}1_{\bbA})=E'$ (see figure \ref{f4}).
This is really just the explicit description of the coreflection of $\bbC$ along the functor ``pullback with $F$''
$$-\times_{\bbB}\bbA\:\Cat(\Q)\slice{\bbB}\to\Cat(\Q).$$
Hence $\Cat(\Q)$ admits all partial products with $F\:\bbA\to\bbB$ if and only if this functor has a right adjoint. S. Niefield [1982] proved that this in turn is equivalent to the functor ``product with $F$'' 
$$-\times F\:\Cat(\Q)\slice{\bbB}\to\Cat(\Q)\slice{\bbB}$$
having a right adjoint, i.e.\ to the {\em exponentiability} of $F$. 
%
%
\par
Suppose now that $F\:\bbA\to\bbB$ and $\bbC$ are given, and that we 
want to construct their partial product $(\bbP,P,E)$. Putting $\bbP'=*_X$ in the 
diagram in figure \ref{f4} and letting $X$ range over all objects 
of $\Q$, the universal property of the partial product dictates at once what the object-set $\bbP_0$ and the 
object-maps $P\:\bbP_0\to\bbC_0$ and 
$E\:(\bbP\times_{\bbB}\bbA)_0\to\bbC_0$ must be:
\begin{itemize}
\item $\bbP_0=\{(b,H)\mid b\in\bbB\mbox{ and }H\:\bbA_b\to\bbC\mbox{ is a functor}\}$, 
with types $t(b,H)=tb$,
\item for $(b,H)\in\bbP_0$, $P(b,H)=b$,
\item for $((b,H),a)\in(\bbP\times_{\bbB}\bbA)_0$, $E((b,H),a)=Ha$.
\end{itemize}
Thus we are left to find a $\Q$-enrichment of the object-set 
$\bbP_0$, making it a $\Q$-category $\bbP$ and making $P$ and $E$ 
functors with the required universal property; the next lemma tells us how to do this.
\begin{lemma}\label{6}
If $F\:\bbA\to\bbB$ satisfies \ref{2}--\ref{a1} and \ref{5}--\ref{a2}, then $\Cat(\Q)$ admits partial products over $F\:\bbA\to\bbB$.
\end{lemma}
\proof
Assuming \ref{2}--\ref{a1} it makes sense to define
\begin{eqnarray*}
\bbP((b',H'),(b,H))
 & \defeq & \bbC(H'-,H-)^F \\
 & = & \mbox{the outcome of applying the right adjoint to the map}\\
 & & \mbox{in \eqref{4.1} on the distributor }\bbC(H'-,H-)\:\bbA_b\dist\bbA_{b'}.
\end{eqnarray*}
Whereas the identity inequality
$$1_{t(b,H)}\leq\bbP((b,H),(b,H))$$
reduces to the fact that $H\:\bbA_b\to\bbC$ is 
a functor, it is the assumed \ref{5}--\ref{a2} together with the composition 
inequality in the $\Q$-category $\bbC$ that assures the composition 
inequality:
$$\bbP((b'',H''),(b',H'))\circ\bbP((b',H'),(b,H))\leq\bbP((b'',H''),(b,H)).$$
This construction clearly makes $P$ and $E$ functorial. As for the universal 
property of $(\bbP,P,E)$, given a $\Q$-category $\bbP'$ and functors 
$P'\:\bbP'\to\bbB$ and $E'\:\bbP'\times_{\bbB}\bbA\to\bbC$, it is 
straightforward to verify that
$$K\:\bbP'\to\bbP\:x\mapsto K(x)\defeq\Big(P'x,\ E'(x,-)\:\bbA_{P'x}\to \bbC\:a\mapsto E'(x,a)\Big)$$
is the required unique factorization.
\endofproof
\par
Finally we shall show that conditions \ref{2}--\ref{a1} and \ref{5}--\ref{a2} are not only sufficient but also necessary for $\Cat(\Q)$ to admit partial products over $F\:\bbA\to\bbB$. Thereto we shall use an auxiliary construction concerning distributors between $\Q$-categories that we better recall beforehand: given a distributor $\Phi\:\bbX\dist\bbY$, we shall say that a co-span of functors like 
$$\xymatrix{\bbX\ar[r]^S & \bbC & \bbY\ar[l]_T}$$
represents $\Phi$ when $\Phi=\bbC(T-,S-)$. Any $\Phi$ admits at least one such representing co-span: let $\bbC_0=\bbX_0\uplus\bbY_0$ and for all $a,a'\in\bbX_0$ and $b,b'\in\bbY_0$ put $\bbC(a',a)=\bbX(a',a)$, $\bbC(b',b)=\bbY(b',b)$, $\bbC(b,a)=\Phi(b,a)$, $\bbC(a,b)=0_{tb,ta}$, so that the co-span of full embeddings
$$\xymatrix{\bbX\ar[r]^{S_{\bbX}} & \bbC & \bbY\ar[l]_{S_{\bbY}}}$$
surely represents $\Phi$. (This latter co-span is universal amongst all representing co-spans for $\Phi$; M. Grandis and R. Par\'e [1999] speak, in the context of double colimits in double categories, of the cotabulator (or gluing, or collage) of $\Phi$. This is however not important for us here; on the contrary, further on it is crucial to consider non-universal representing co-spans.)
\begin{lemma}\label{6.1}
If $\Cat(\Q)$ admits partial products over $F\:\bbA\to\bbB$, then $F\:\bbA\to\bbB$ satisfies \ref{2}--\ref{a1} and \ref{2}--\ref{a2}.
\end{lemma}
\proof
For $b,b'\in\bbB$ and $\Phi\:\bbA_b\dist\bbA_{b'}$, choose a representing co-span 
$$\xymatrix{\bbA_b\ar[r]^S & \bbC &\bbA_{b'}\ar[l]_T}.$$
Considering the partial product of $F$ with $\bbC$, say $(\bbP,P,E)$, it is a fact that the hom-arrow $\bbP((b',T),(b,S))$ is a $\Q$-arrow smaller than $\bbB(b',b)$. Now, any $\Q$-arrow $f\:X\to Y$ determines a $\Q$-category\footnote{This is actually an instance of the universal representing co-span, when viewing the $\Q$-arrow $f\:X\to Y$ as a one-element distributor $(f)\:*_X\dist *_Y$ between one-object $\Q$-categories.} $\bbP_f$ like so:
\begin{itemize}
\item objects: $(\bbP_f)_0=\{X\}\uplus\{Y\}$ with $tX=X\in\Q$ and $tY=Y\in\Q$,
\item hom-arrows: $\bbP_f(Y,X)=f$, $\bbP_f(X,X)=1_X$, $\bbP_f(Y,Y)=1_Y$ and $\bbP_f(X,Y)=0_{Y,X}$.
\end{itemize}
The inequality $f\leq\bbB(b',b)$ holds if and only if 
$$P_f\:\bbP_f\to\bbB\:X\mapsto b,Y\mapsto b'$$ 
is a functor; and similarly the collection of inequalities $f\wedge\bbA(a',a)\leq\Phi(a',a)$ (one for each $a\in\bbA_b,a'\in\bbA_{b'}$) 
is equivalent to 
$$E_f\:\bbP_f\times_{\bbB}\bbA\to\bbC\:(X,a)\mapsto a,(Y,a')\mapsto a'$$ 
being a functor. Using the universal property of the partial product $(\bbP,P,E)$ one easily checks that $P_f$ and $E_f$ 
determine and are determined by the single functor 
$$K\:\bbP_f\to\bbP\:X\mapsto(b,S),Y\mapsto(b',T),$$
whose functoriality in turn is equivalent to the inequality $f\leq\bbP((b',T),(b,S))$.
\par
The above argument is actually independent of the chosen representing co-span for $\Phi$: if another co-span
$$\xymatrix{\bbA_b\ar[r]^{S'} & \bbC' &\bbA_{b'}\ar[l]_{T'}}$$
also represents $\Phi$, and $(\bbP',P',E')$ denotes the partial product of $F$ with $\bbC'$, then the ``same'' argument shows that, for any $\Q$-arrow $f\leq\bbB(b',b)$, the collection of inequalities $f\wedge\bbA(a',a)\leq\Phi(a',a)$ (one for each $a\in\bbA_b,a'\in\bbA_{b'}$) 
is equivalent to the single inequality $f\leq\bbP'((b',T'),(b,S'))$. Thus it follows that $\bbP((b',T),(b,S))=\bbP'((b',T'),(b,S'))$.
\par
As a result the map
\begin{equation}\label{6.0}
\Dist(\Q)(\bbA_{b},\bbA_{b'})\to\downarrow\bbB(b',b)\:\Phi\mapsto\Phi^F\defeq\bbP((b',T),(b,S)),
\end{equation}
where one computes $\Phi^F$ with the aid of {\it any chosen representing co-span} for $\Phi$, is (well-defined and) the right adjoint in \eqref{4.1b}. We end by showing that it satisfies the lax commutativity of the diagram in figure \ref{f3}; thereto it is important that, in the map prescription of \eqref{6.0}, any chosen representing co-span for a given distributor will do.
\par
For $b,b',b''\in\bbB$, $\Phi\:\bbA_b\dist\bbA_{b'}$ and 
$\Psi\:\bbA_{b'}\dist\bbA_{b''}$, consider the $\Q$-category $\bbC$ like so:
\begin{itemize}
\item objects: $\bbC_0=(\bbA_b)_0\uplus(\bbA_{b'})_0\uplus(\bbA_{b''})_0$ with ``inherited types'',
\item hom-arrows: for all $a,a_1\in\bbA_b$, $a',a_1'\in\bbA_{b'}$ and $a'',a_1''\in\bbA_{b''}$, put $\bbC(a_1,a)=\bbA_{b}(a_1,a)$, $\bbC(a_1',a')=\bbA_{b'}(a_1',a')$, $\bbC(a_1'',a'')=\bbA_{b''}(a_1'',a'')$, $\bbC(a',a)=\Phi(a',a)$, $\bbC(a'',a')=\Psi(a'',a')$ and $\bbC(a'',a)=(\Psi\tensor\Phi)(a'',a)$, all other 
hom-arrows are zero.
\end{itemize}
The co-spans of full embeddings
$$\xymatrix{\bbA_b\ar[r]^{S} & \bbC & \bbA_{b'}\ar[l]_T}, 
\xymatrix{\bbA_{b'}\ar[r]^{T} & \bbC & \bbA_{b''}\ar[l]_U},
\xymatrix{\bbA_{b}\ar[r]^{S} & \bbC & \bbA_{b''}\ar[l]_U}$$
represent respectively $\Phi$, $\Psi$ and $\Psi\tensor\Phi$. Writing $(\bbP,P,E)$ for the partial product of $F$ and $\bbC$, the  
compostion-inequality
$$\bbP((b'',U),(b',T))\circ\bbP((b',T),(b,S))\leq\bbP((b'',U),(b,S))$$
says precisely that $\Psi^F\circ\Phi^F\leq(\Psi\tensor\Phi)^F$, as wanted.
\endofproof




\end{document}